\documentclass[12pt]{amsproc}
  \usepackage[margin=1in]{geometry}
\usepackage{setspace,fullpage}
\usepackage{graphicx}
\usepackage[nice]{nicefrac}
\usepackage{amssymb}
\usepackage{amsmath,amsthm,amscd,amssymb, mathrsfs}

\usepackage{latexsym}
\usepackage[colorlinks,citecolor=red,pagebackref,hypertexnames=false]{hyperref}
\usepackage{graphicx}
\usepackage{subfigure}
\usepackage{CJKutf8}

\usepackage{graphicx}
\usepackage{hyperref}
\usepackage{cite}
\usepackage{indentfirst}

\newcommand\inn[2]{\langle#1,#2 \rangle}
\newcommand\linn[2]{\left\langle#1,#2 \right\rangle}
\newcommand\innerHaar[2]{\widehat{#1}\left( #2 \right)}
\newcommand\inner[1]{\langle#1 \rangle}
  \newcommand{\norm}[1]{\left\lVert#1\right\rVert}
\newcommand\bracket[1]{\left(#1\right)}
\newcommand\cbracket[1]{\left[#1\right]}
\newcommand{\Atwo}{[\omega]_{A_2}}
\newcommand{\Asum}{\sum_{I\in\mathcal{D}}\sum_\alpha}
\newcommand{\Bsum}{\sum_{J\in\mathcal{D}}\sum_\beta}
\newcommand{\wHaara}{h_{E_{\alpha,I}}}

\newcommand{\Ea}{{E_{\alpha,I}}}
\newcommand{\Eb}{{E_{\beta,J}}}
\newcommand{\Carlsum}{\sup_{\alpha,I}\frac{1}{|\Ea|}\sum_{J\subset I}\sum_{\beta:\Eb\subset\Ea}}
\newcommand{\wCarlsum}{\sup_{\alpha,I}\frac{1}{|\Ea\inner{\omega}_\Ea|}\sum_{J\subset I}\sum_{\beta:\Eb\subset\Ea}}
\newcommand{\wCarlsumne}{\sup_{\alpha,I}\frac{1}{|\Ea\inner{\omega^{-1}}_\Ea|}\sum_{J\subset I}\sum_{\beta:\Eb\subset\Ea}}
\newtheorem{thm}{Theorem}[section]
\newtheorem{defn}[thm]{Definition}

\newtheorem{lem}[thm]{Lemma}
\newtheorem{prop}[thm]{Proposition}

\title{A note on sharp weighted bound for Haar shift and multiplier}
\author[]{Chih-Chieh Hung}
\address{Department of Mathematics\\ National Taiwan University\\
Taipei, 106 Taiwan}
\email{r07221003@ntu.edu.tw}

\author[C-Y. Shen]{Chun-Yen Shen}
\address{Department of Mathematics\\ National Taiwan University and National Center for Theoretical Sciences\\Taipei, 106 Taiwan}
\email{cyshen@math.ntu.edu.tw}

\thanks{supported by MOST through grant 108-2628-M-002-010-MY4}

\begin{document}
\begin{abstract}
    We provide  elementary proofs for the terms that are left in the work of Kelly Bickel, Sandra Pott, Maria Carmen Reguera, Eric T. Sawyer and Brett D. Wick that proved the sharp weighted $A_2$ bound for Haar shifts and Haar multiplier paraproducts. Our proofs use weighted square function estimate, Carleson embedding and Wilson's system. 
\end{abstract}

\maketitle

\section{Introduction}
One of the recent important results in the area of Harmonic analysis is the solution of the longstanding $A_2$ conjecture concerning with the sharp weighted bound for Calder\'on-Zygmund singular integrals. There are many people who have made great efforts to the $A_2$ problem including S. Buckley, S. Petermichl, A. Volberg, F. Nazarov, S. Treil, to name just a few. Finally T. Hyt\"{o}nen proved the $A_2$ conjecture for general Calder\'on-Zygmund singular integrals. After that, an interesting question is to find some elementary proofs for the sharp $A_2$ bound . It was  Lerner who first found a simple proof to the $A_2$ conjecture by 
using Sparse operator\cite{Lerner12}. Yet another simple and different proof for Hilbert transform \cite{Haar} was given by S. Pott, M. Reguera, E. Sawyer and B. Wick. Their proof reduced proving the sharp $A_2$ bound for the Hilbert transform to estimate nine terms which are called canonical individual paraproduct composition operators. More precisely, the linear bound of Hilbert transform in the $A_2$ characteristic:
\begin{equation*}
    \norm{H}_{L^2(\omega)\rightarrow L^2(\omega)}\lesssim\Atwo
\end{equation*}
can be reduced to the nine component operators in the natural weighted resolution of the conjugation $M_{\omega^\frac{1}{2}}SM_{\omega^\frac{-1}{2}}$ induced by a multiplier into paraproducts where $S$ is a Haar shift operator and 

\begin{equation*}
    M_{\omega^\frac{1}{2}}SM_{\omega^\frac{-1}{2}}=\bracket{P^{(0,1)}_{\widehat{\omega^{\frac{1}{2}}}}+P^{(1,0)}_{\widehat{\omega^{\frac{1}{2}}}}+P^{(0,0)}_{\inner{\omega^\frac{1}{2}}}}S\bracket{P^{(0,1)}_{\widehat{\omega^{\frac{-1}{2}}}}+P^{(1,0)}_{\widehat{\omega^{\frac{-1}{2}}}}+P^{(0,0)}_{\inner{\omega^\frac{-1}{2}}}}
\end{equation*}
(see section \ref{PreHaar} for the definition of these operators).
In their work, eight of the nine operators above $P^{(\alpha,\beta)}_{\widehat{\omega^{\frac{1}{2}}}}SP^{(\gamma,\delta)}_{\widehat{\omega^{\frac{-1}{2}}}}$ they were able to use some elementary proofs to show their norms are bounded by $\Atwo$. However an operator was left out without giving elementary proofs but directly appealing the already known bound by S. Petermichl's result\cite{Hsharp} for the Hilbert transform. Hence they left an open question whether one can use elementary proofs to give the sharp bound for this operator.

Now for the Haar multiplier, K. Bickel, E. Sawyer and B. Wick further developed these elementary proofs for the Haar multiplier in $\mathbb{R}^d$ (see section \ref{PreMul} for more details). They also decomposed the Haar multiplier into nine component operators. However, just like the situation in Haar shifts above, there was an operator left out without giving elementary proofs but directly appealing the already known bound by D. Chung.
\cite{Chung}. Therefore they left an open question whether one can use elementary proofs to give the sharp bound for this operator.

The purpose of this paper is to give elementary proofs for these two questions. In other words, we prove
\begin{equation*}
    \norm{P^{(0,1)}_{\widehat{\omega^{\frac{1}{2}}}}SP^{(1,0)}_{\widehat{\omega^{\frac{-1}{2}}}}\phi}_{L^2}\lesssim \Atwo,
\end{equation*}
and
\begin{equation*}
    \norm{P^{(0,1)}_{\widehat{\omega^{\frac{1}{2}}}}T_\sigma P^{(1,0)}_{\widehat{\omega^{\frac{-1}{2}}}}\phi}\lesssim \Atwo,
\end{equation*}
which answer the questions in \cite{Haar} and \cite{Multiplier}.
 We will define the notations for Haar shifts and Haar multiplier in section \ref{PreHaar} and \ref{PreMul} respectively. Meanwhile in order to keep our proof as simple as possible, our strategy is to combine three operators of the nine operators and redecompose into new operators, and each of them we are able to prove the sharp $A_2$ bound.

\section{Preliminary for Haar shift}\label{PreHaar}
Let $L^2\equiv L^2(\mathbb{R})$ be the space of square integrable functions over $\mathbb{R}$. A weight is a nonnegative local integrable function i.e. $\omega\in L^1_{loc}(\mathbb{R})$ and $$\omega:\mathbb{R}\rightarrow \mathbb{R}_{\geq 0}.$$
Also $\inn{}{}:=\inn{}{}_{L^2(\mathbb{R})}$, and $\inn{}{}_\omega:=\inn{}{}_{L^2(\mathbb{R},\omega)}$; $\norm{\cdot}=\norm{\cdot}_{L^2(\mathbb{R})}$, and $\norm{\cdot}_\omega=\norm{\cdot}_{L^2(\mathbb{R},\omega)}$. $\mathcal{D}$ is the usual dyadic grid of intervals on the real line, and $S$ is Haar shifts defined by $Sh_I=h_{I-}\quad\forall I\in\mathcal{D}$. And the $A_2$ characteristic of $\omega$ is defined by:
$$\Atwo\equiv\sup_{I\in\mathcal{D}}\inner{\omega}_I\inner{\omega^{-1}}_I,$$
where $\inner{\omega}_I$ denotes the average of $\omega$ over a cube $I\in\mathcal{D}$. 

\subsection{Haar function}
Define the Haar function $h^0_I$ and the average function $h^1_I$ by 
\begin{equation*}
    h^0_I\equiv h_I\equiv\frac{1}{\sqrt{I}}\bracket{1_{I-}-1_{I+}}\quad\textit{ and }h^1_I\equiv\frac{1}{|I|}1_I,\quad I\in\mathcal{D}.
\end{equation*}
The paraproduct operators considered in this part are the following dyadic operators which come from \cite{Haar}.
\begin{defn}
Given a  symbol $b=\{b_I\}_{I\in\mathcal{D}}$ and a pair $(\alpha,\beta)\in\{0,1\}\times\{0,1\}$, define dyadic paraproduct acting on a function $f$ by
\begin{equation*}
    P^{\alpha,\beta}_bf\equiv\sum_{I\in\mathcal{D}}b_I\inn{f}{h^\beta_I}h^\alpha_I.
\end{equation*}
\end{defn}
For a function $b$ and $I\in\mathcal{D}$ we let 
\begin{align*}
    &\innerHaar{b}{I}\equiv\inn{b}{h^0_I}\\
    &\inner{b}_I\equiv\inn{b}{h^1_I},
\end{align*}
One of the key tools we need is disbalanced Haar functions. Given a weight $\omega$ on $\mathbb{R}$ we let
\begin{equation*}
    C_{K}(\omega):=\sqrt{\frac{\inner{\omega}_{K+}\inner{\omega}_{K-}}{\inner{\omega}_K}}\quad D_{K}(\omega):=\frac{\innerHaar{\omega}{K}}{\inner{\omega}_K}.
\end{equation*}
Then we have
\begin{align}\label{Haardec}
    &\inn{h_{K}}{g}_\omega\\
    =&\sqrt{\frac{\inner{\omega}_{K+}\inner{\omega}_{K-}}{\inner{\omega}_K}}\inn{h^\omega_{K}}{g}_\omega+\frac{\innerHaar{\omega}{K}}{\inner{\omega}_K}\inn{h^1_{K}}{g}_\omega\\
    =&C_{K}(\omega)\inn{h^\omega_{K}}{g}_\omega+D_{K}(\omega)\inn{h^1_{K}}{g}_\omega.
\end{align}

Also we have estimates for $C_K(\omega)$. By definition, we have
\begin{equation}\label{Cestimate}
    C_K(\omega)\lesssim \sqrt{\inner{\omega}_K}
\end{equation}
and
\begin{equation}\label{C-estimate}
    C_{K-}(\omega)\lesssim \sqrt{\inner{\omega}_K}.
\end{equation}

\subsection{Carleson Embedding}
The following is the classical dyadic Carleson embedding theorem (for example see \cite{Pereyra}).
Let $\lambda=\{\lambda_I:I\in\mathcal{D}\}$ be a sequence of positive numbers. Define 
$$\norm{\lambda}_{CM}:=\sup_{J}\frac{1}{|J|}\sum_{I\in\mathcal{D}(J)}\lambda_I,$$

and we say that a sequence $\{\lambda_I\}_{I\in\mathscr{D}}$ of positive numbers is a Carleson sequence if there exists a constant $C>0$ such that
$$\norm{\lambda}_{CM}\leq C.$$
\begin{thm}\label{Car_emb}(Carleson's embedding theorem) Given a Carleson sequence $\{\lambda_I\}_{I\in\mathscr{D}}$, then for all $f\in L^2(\mathbb{R})$
$$\Sigma_{I\in\mathscr{D}}|\inner{f}_I|^2\lambda_I\leq C\norm{f}^2_2.$$
\end{thm}

We also use the following weighted version of Carleson embedding theorem (for example see the work of F. Nazarov, S. Treil, and A. Volberg\cite{NTV97}).
\begin{thm}\label{wCarleson}
(Weighted Carleson embedding)\\
The following are equivalent
\begin{align}
    \sum_I\inner{f\omega^{1/2}}^2_I\alpha_I\leq C\norm{f}^2\\
    \forall J\in\mathcal{D},\quad \frac{1}{|J|}\sum_{I\subset J}\inner{\omega}^2_I\alpha_I\leq C\inner{\omega}_J.
\end{align}
\end{thm}
Let $f=g\omega^{1/2}$, we get
\begin{thm}\label{Carlesonem}
The following are equivalent
\begin{align}
    \sum_I\alpha_I\inner{g\omega}^2_I=\sum_I\alpha_I\inn{g}{h^1_I}_\omega^2\leq C\norm{g}_\omega^2\\
    \forall J\in\mathcal{D},\quad \frac{1}{|J|}\sum_{I\subset J}\alpha_I\inner{\omega}^2_I\leq C\inner{\omega}_J.
\end{align}
\end{thm}
We also have the following bilinear embedding result \cite{Hsharp}.
\begin{thm}\label{Bilinear Imbedding Theorem}(Bilinear Embedding Theorem)\\
Let $\omega$ and $\nu$ be weights so that $\inner{\omega}_I\inner{\nu}_I\leq C$ for all intervals $I$ and let $\alpha_I$ be nonnegative sequence so that the three estimates below hold for all $J$
\begin{align}
    &\sum_{I\subset J}\alpha_I\inner{\nu}_I\leq C\nu(J)\\
    &\sum_{I\subset J}\alpha_I\inner{\omega}_I\leq C\omega(J)\\
    &\sum_{I\subset J}\alpha_I\leq C|J|.
\end{align}
Then for all $f\in L^2(\omega)$ and $g\in L^2(\nu)$
$$\sum_I\alpha_I\inn{f}{h^1_I}_\omega\inn{g}{h^1_I}_\nu\lesssim C \norm{f}_\omega\norm{g}_\nu.$$
\end{thm}
Also we have the following inequalities from S. Petermichl's work \cite{Hsharp}
\begin{align}
    &\sum_{I\subset J}\frac{|\innerHaar{\omega^{-1}}{I}\innerHaar{\omega}{I-}|}{\inner{\omega}_I}\lesssim {\Atwo}\omega^{-1}(J)\label{wineq}\\
    &\sum_{I\subset J}\frac{|\innerHaar{\omega^{-1}}{I}\innerHaar{\omega}{I-}|}{\inner{\omega^{-1}}_I}\lesssim {\Atwo}\omega(J)\label{iwineq}\\
    &\sum_{I\subset J}|\innerHaar{\omega^{-1}}{I}\innerHaar{\omega}{I-}|\lesssim {\Atwo}|J|\label{ineq}.
\end{align}
\subsection{Square function}
The dyadic square function is given by
\begin{equation*}
    Sf(x)\equiv\sqrt{\sum_{I\in\mathcal{D}}|\innerHaar{f}{I}|^2h^1_I(x)},
\end{equation*}
and for any weight $\sigma\geq 0$ we have
\begin{equation*}
    \norm{Sf}^2_\sigma=\sum_{I\in\mathcal{D}}|\innerHaar{f}{I}|^2\inner{\sigma}_I.
\end{equation*}
It has been shown in \cite{Hsquare} that if $\omega\in A_2$, then
\begin{equation}\label{Swestimate}
    \norm{Sf}_\omega\lesssim\Atwo\norm{f}_\omega.
\end{equation}
Applying this equality to $f=\omega^\frac{-1}{2}1_I$ for $I\in\mathcal{D}$ will lead to
\begin{align*}
    \sum_{J\subset I}\innerHaar{\omega^\frac{-1}{2}}{J}^2\inner{\omega}_J\lesssim\Atwo^2|I|\quad\forall I\in\mathcal{D}.
\end{align*}
There is also a modified version of square function which incorporates a shift in the indices from \cite{Haar}. Define the modified square function $S_\pi$ by,
\begin{equation*}
    S_\pi f(x)\equiv\sqrt{\sum_{I\in\mathcal{D}}|\innerHaar{f}{I}|^2\frac{1}{|I|}1_{\pi I}(x)},
\end{equation*}
where $\pi I$ is the dyadic parent of the dyadic interval $I.$
\begin{thm}\label{modified}\cite{Haar}
For any $f\in L^2(\omega)$, we have
\begin{equation}
    \norm{S_\pi \phi}_\omega\lesssim\Atwo\norm{\phi}_\omega.
\end{equation}
\end{thm}
\section{Main result for Haar shift}
What we will prove is the following 

$$\norm{P^{(0,1)}_{\widehat{\omega^{\frac{1}{2}}}}SP^{(1,0)}_{\widehat{\omega^{\frac{-1}{2}}}}\phi}\lesssim \Atwo\norm{\phi},$$

it is equivalent to show

$$\norm{\omega^\frac{1}{2} SP^{(1,0)}_{\widehat{\omega^{\frac{-1}{2}}}}\phi}\lesssim \Atwo\norm{\phi}$$
since  
\begin{align*}
    \bracket{P^{(0,1)}_{\widehat{\omega^{\frac{1}{2}}}}+P^{(1,0)}_{\widehat{\omega^{\frac{1}{2}}}}+P^{(0,0)}_{\inner{\omega^\frac{1}{2}}}}SP^{(1,0)}_{\widehat{\omega^{\frac{-1}{2}}}}\phi=\omega^\frac{1}{2} SP^{(1,0)}_{\widehat{\omega^{\frac{-1}{2}}}}\phi
\end{align*}
and it has been shown in \cite{Haar} that we already have
\begin{align*}
    \norm{P^{(1,0)}_{\widehat{\omega^{\frac{1}{2}}}}SP^{(1,0)}_{\widehat{\omega^{\frac{-1}{2}}}}}\lesssim\Atwo\\
    \norm{P^{(0,0)}_{\inner{\omega^\frac{1}{2}}}SP^{(1,0)}_{\widehat{\omega^{\frac{-1}{2}}}}}\lesssim\Atwo.
\end{align*}
Also by duality we are going to prove that $\forall g\in L^2(\omega)$
$$\linn{SP^{(1,0)}_{\widehat{\omega^{\frac{-1}{2}}}}\phi}{g}_\omega\lesssim\Atwo\norm{\phi}\norm{g}_\omega.$$

By definition,
$$P^{(1,0)}_{\widehat{\omega^{\frac{-1}{2}}}}\phi=\sum_J\innerHaar{\omega^{\frac{-1}{2}}}{J}\innerHaar{\phi}{J}h_J$$

and 
\begin{align*}
    SP^{(1,0)}_{\widehat{\omega^{\frac{-1}{2}}}}\phi&=\sum_K\bracket{\sum_{J\subsetneq K}\innerHaar{\omega^{\frac{-1}{2}}}{J}\innerHaar{\phi}{J}\inn{h^1_J}{h_K}}h_{K-}\\
    &=\sum_K\bracket{\innerHaar{\omega^{\frac{-1}{2}}\phi}{K}-\inner{\omega^{\frac{-1}{2}}}_K\innerHaar{\phi}{K}-\inner{\phi}_K\innerHaar{\omega^{\frac{-1}{2}}}{K}}h_{K-}
\end{align*}
\begin{align*}
    &\inn{SP^{(1,0)}_{\widehat{\omega^{\frac{-1}{2}}}}\phi}{g}_\omega\\
    =&\sum_K\bracket{\innerHaar{\omega^{\frac{-1}{2}}\phi}{K}-\inner{\omega^{\frac{-1}{2}}}_K\innerHaar{\phi}{K}-\inner{\phi}_K\innerHaar{\omega^{\frac{-1}{2}}}{K}}\inn{h_{K-}}{g}_\omega\\
    =&\sum_K\bracket{\innerHaar{\omega^{\frac{-1}{2}}\phi}{K}-\inner{\omega^{\frac{-1}{2}}}_K\innerHaar{\phi}{K}-\inner{\phi}_K\innerHaar{\omega^{\frac{-1}{2}}}{K}}\bracket{C_{K-}(\omega)\inn{h^\omega_{K-}}{g}_\omega+D_{K-}(\omega)\inn{h^1_{K-}}{g}_\omega}\\
    =&\sum_K\innerHaar{\omega^{\frac{-1}{2}}\phi}{K}C_{K-}(\omega)\inn{h^\omega_{K-}}{g}_\omega\tag{A1}\\
    -&\sum_K\inner{\omega^{\frac{-1}{2}}}_K\innerHaar{\phi}{K}C_{K-}(\omega)\inn{h^\omega_{K-}}{g}_\omega\tag{B1}\\
    -&\sum_K\inner{\phi}_K\innerHaar{\omega^{\frac{-1}{2}}}{K}C_{K-}(\omega)\inn{h^\omega_{K-}}{g}_\omega\tag{C1}\\
    +&\sum_K\innerHaar{\omega^{\frac{-1}{2}}\phi}{K}D_{K-}(\omega)\inn{h^1_{K-}}{g}_\omega\tag{A2}\\
    -&\sum_K\inner{\omega^{\frac{-1}{2}}}_K\innerHaar{\phi}{K}D_{K-}(\omega)\inn{h^1_{K-}}{g}_\omega\tag{B2}\\
    -&\sum_K\inner{\phi}_K\innerHaar{\omega^{\frac{-1}{2}}}{K}D_{K-}(\omega)\inn{h^1_{K-}}{g}_\omega.\tag{C2}
\end{align*}
Therefore 
\begin{equation*}
    \norm{SP^{(1,0)}_{\widehat{\omega^{\frac{-1}{2}}}}\phi}_\omega\lesssim \Atwo\norm{\phi}
\end{equation*}
is equivalent to show that these six terms A1, A2, B1, B2, C1, C2 are all bounded by $\Atwo\norm{\phi}\norm{g}_\omega$.
\paragraph{For B1 part}
\begin{align*}
    &\sum_K\inner{\omega^{\frac{-1}{2}}}_K\innerHaar{\phi}{K}C_{K-}(\omega)\inn{h^\omega_{K-}}{g}_\omega\\
    \leq& \bracket{\sum_K\inner{\omega^{\frac{-1}{2}}}_K\innerHaar{\phi}{K}^2|C_{K-}(\omega)|}^{\frac{1}{2}}\bracket{\sum_K\inner{\omega^{\frac{-1}{2}}}_K|C_{K-}(\omega)|\inn{h^\omega_{K-}}{g}^2_\omega}^{\frac{1}{2}}\quad \textit{by Cauchy-Schwarz}\\
    \lesssim& \bracket{\sum_K\sqrt{\inner{\omega}_K\inner{\omega^{-1}_K}}\innerHaar{\phi}{K}^2}^{\frac{1}{2}}\bracket{\sum_K\sqrt{\inner{\omega}_K\inner{\omega^{-1}_K}}\inn{h^\omega_{K-}}{g}^2_\omega}^{\frac{1}{2}}\quad \textit{by (1) and (2)}\\
    \leq& \Atwo^{1/2}\norm{\phi}\norm{g}_\omega\quad \textit{by definition of $A_2$.}
\end{align*}
\paragraph{For C1 part} 
\begin{align*}
    &\sum_K\inner{\phi}_K\innerHaar{\omega^{\frac{-1}{2}}}{K}C_{K-}(\omega)\inn{h^\omega_{K-}}{g}_\omega\\
    \leq&\bracket{\sum_K\innerHaar{\omega^{\frac{-1}{2}}}{K}^2C_{K-}(\omega)^2\inner{\phi}_K^2}^{1/2}\norm{g}_\omega \quad\textit{by Cauchy-Schwarz}\\
    \lesssim&\bracket{\sum_K\innerHaar{\omega^{\frac{-1}{2}}}{K}^2\inner{\omega}_K\inner{\phi}_K^2}^{1/2}\norm{g}_\omega \quad\textit{by (2)}\\
    \lesssim&\norm{g}_\omega\norm{\phi}\bracket{\sup\frac{1}{|J|}\sum_{K\subset J}\innerHaar{\omega^{\frac{-1}{2}}}{K}^2\inner{\omega}_K}^{1/2}\quad\textit{by Carleson Embedding}\\
    \lesssim&\Atwo\norm{g}_\omega\norm{\phi}\quad\textit{by sharp bound of square function.}
\end{align*}
\paragraph{For B2 part}
\begin{align*}
    &\sum_K\inner{\omega^{\frac{-1}{2}}}_K\innerHaar{\phi}{K}D_{K-}(\omega)\inn{h^1_{K-}}{g}_\omega\\
    &\leq \norm{\phi}\bracket{\sum_K\inner{\omega^{\frac{-1}{2}}}^2_K\bracket{\frac{\innerHaar{\omega}{K-}}{\inner{\omega}_{K-}}}^2\inner{g\omega}_{K-}^2}^\frac{1}{2}\quad\textit{by Cauchy-Schwarz}\\
    &\lesssim \norm{\phi}\norm{g}_\omega\bracket{\sup_J \frac{1}{\inner{\omega}_J|J|}\sum_{I\subset J} \inner{\omega^{\frac{-1}{2}}}^2\innerHaar{\omega}{K-}^2}^\frac{1}{2},\quad\textit{by Theorem \ref{Carlesonem}}
\end{align*}
where
\begin{align*}
    &\sup_J \frac{1}{\inner{\omega}_J|J|}\sum_{I\subset J} \inner{\omega^{\frac{-1}{2}}}^2\innerHaar{\omega}{K-}^2\\
    &\lesssim\Atwo^2\sup_J \frac{1}{\inner{\omega}_J|J|}\int_J \omega\quad\textit{by theorem \ref{modified}.}\\
    &\lesssim\Atwo.
\end{align*}
So this term is bounded by
$$\Atwo\norm{\phi}\norm{g}_\omega.$$
\paragraph{For C2 part}
In this part, we need to split it into two cases. Since we have
\begin{align*}
    \innerHaar{\omega^{\frac{-1}{2}}}{K}=C_K(\omega^{-1})\inn{\omega^{1/2}}{h^{\omega^{-1}}_K}_{\omega^{-1}}+D_K(\omega^{-1})\inner{\omega^{-1/2}}_K,\quad\textit{by (\ref{Haardec})}
\end{align*}
we get the following 
\begin{align*}
    &\sum_K\inner{\phi}_K\innerHaar{\omega^{\frac{-1}{2}}}{K}D_{K-}(\omega)\inn{h^1_{K-}}{g}_\omega\tag{C2}\\
    =&\sum_KC_K(\omega^{-1})\inn{\omega^{1/2}}{h^{\omega^{-1}}_K}_{\omega^{-1}}\inner{\phi}_KD_{K-}(\omega)\inn{h^1_{K-}}{g}_\omega\tag{C21}\\
    +&\sum_KD_K(\omega^{-1})\inner{\omega^{-1/2}}_K\inner{\phi}_KD_{K-}(\omega)\inn{h^1_{K-}}{g}_\omega.\tag{C22}
\end{align*}
\paragraph{For C21}
\begin{align*}
    &\sum_KC_K(\omega^{-1})\inn{\omega^{1/2}}{h^{\omega^{-1}}_K}_{\omega^{-1}}\inner{\phi}_KD_{K-}(\omega)\inn{h^1_{K-}}{g}_\omega\\
    \leq&\bracket{\sum_K\inn{\omega^{1/2}}{h_K^{\omega^{-1}}}^2_{\omega^{-1}}\inner{\phi}^2_K}^{1/2}\bracket{\sum_K\inner{\omega^{-1}}_K\frac{\innerHaar{\omega}{K-}^2}{\inner{\omega}_{K-}^2}\inn{h^1_{K-}}{g}_\omega^2}^{1/2}\textit{by Cauchy-Schwarz}\\ 
    \leq &\bracket{\sup_J \frac{1}{|J|}\sum_{K\subset J}\inn{\omega^{1/2}}{h_K^{\omega^{-1}}}^2_{\omega^{-1}}}^{1/2}\norm{\phi}\bracket{\sup_J\frac{1}{|J|\inner{\omega}_J}\sum_{K\subset J}\inner{\omega^{-1}}_K\innerHaar{\omega}{K-}^2}^{1/2}\norm{g}_\omega.
\end{align*}
The last inequality comes from theorem (\ref{Car_emb}) and (\ref{Carlesonem}). Note that
\begin{align*}
    &\sup_J \frac{1}{|J|}\sum_{K\subset J}\inn{\omega^{1/2}}{h_K^{\omega^{-1}}}^2_{\omega^{-1}}\\
    \leq &\sup_J \frac{1}{|J|}\sum_{K\cap J\neq\emptyset}\inn{\omega^{1/2}}{h_K^{\omega^{-1}}}^2_{\omega^{-1}}\\
    =&\sup_J\frac{1}{|J|}\norm{\omega^{1/2}\chi_J}_{\omega^{-1}}^2\\
    =&1.
\end{align*}
Also by weighted bound of modified square function (\ref{modified}), we have
\begin{align*}
    &\sup_J\frac{1}{|J|\inner{\omega}_J}\sum_K\inner{\omega^{-1}}_K\innerHaar{\omega}{K-}^2\\
    \leq &\sup_J\frac{1}{|J|\inner{\omega}_J}\sum_{K\cap J\neq\emptyset}\inner{\omega^{-1}}_K\innerHaar{\omega}{K-}^2\\
    \leq&\sup_J \frac{1}{|J|\inner{\omega}_J}\Atwo\norm{\omega\chi_J}_{\omega^{-1}}^2\\
    =&\Atwo.
\end{align*}
So
\begin{align*}
    &\sum_KC_K(\omega^{-1})\inn{\omega^{1/2}}{h^{\omega^{-1}}_K}_{\omega^{-1}}\inner{\phi}_KD_{K-}(\omega)\inn{h^1_{K-}}{g}_\omega\\
    \lesssim& \Atwo\norm{\phi}\norm{g}_\omega.
\end{align*}
\paragraph{For C22}
By Cauchy-Schwarz,
\begin{align*}
    &\sum_KD_K(\omega^{-1})\inner{\omega^{-1/2}}_K\inner{\phi}_KD_{K-}(\omega)\inn{h^1_{K-}}{g}_\omega\\
    \leq &\bracket{\sum_K\frac{|\innerHaar{\omega^{-1}}{K}\innerHaar{\omega}{K-}|}{\inner{\omega^{-1}}_K}\inner{\omega^{-1/2}}^2_K\inner{\phi}^2}^{1/2}\bracket{\sum_K\frac{\innerHaar{\omega^{-1}}{K}\innerHaar{\omega}{K-}}{\inner{\omega^{-1}}_K\inner{\omega}^2_{K-}}\inn{g}{h^1_{K-}}_\omega^2}^{1/2}\\
    \leq& \bracket{\sup \frac{1}{|J|}\sum_{K\subset J}\frac{|\innerHaar{\omega^{-1}}{K}\innerHaar{\omega}{K-}|}{\inner{\omega^{-1}}_K}\inner{\omega^{-1/2}}^2_K}^{1/2}\norm{\phi}\bracket{\sup_J\frac{1}{|J|}\sum_{K\subset J}\frac{\innerHaar{\omega^{-1}}{K}\innerHaar{\omega}{K-}}{\inner{\omega^{-1}}_K}}^{1/2}\norm{g}_\omega\quad\textit{(\ref{Car_emb}), (\ref{Carlesonem})}\\
    \leq& \bracket{\sup \frac{1}{|J|}\sum_{K\subset J}\frac{|\innerHaar{\omega^{-1}}{K}\innerHaar{\omega}{K-}|}{\inner{\omega^{-1}}_K}\inner{\omega^{-1/2}}^2_K}^{1/2}\norm{\phi}\Atwo\norm{g}_\omega\quad\textit{by (\ref{wineq})}\\
    \lesssim& \bracket{\sup \frac{1}{|J|}\sum_{K\subset J}|\innerHaar{\omega^{-1}}{K}\innerHaar{\omega}{K-}|}^{1/2}\norm{\phi}[\omega]^{1/2}_{\Atwo}\norm{g}_\omega\\
    \lesssim&\Atwo\norm{\phi}\norm{g}_\omega\quad\textit{by (\ref{ineq})}.
\end{align*}
\paragraph{For A part}
Now we split $A_1$, and $A_2$ into four parts. Note that 
\begin{align*}
    &\inn{\omega^{-1/2}\phi}{h_K}\\
    &=C_K(\omega^{-1})\inn{\omega^{-1/2}\phi}{h^{\omega^{-1}}_K}+D_K(\omega^{-1})\inner{\phi\omega^{-1/2}}_K\quad\textit{by (\ref{Haardec})}.
\end{align*}
Thus 
\begin{align*}
    &\sum_K\innerHaar{\omega^{\frac{-1}{2}}\phi}{K}C_{K-}(\omega)\inn{h^\omega_{K-}}{g}_\omega\tag{A1}\\
    =&\sum_KC_K(\omega^{-1})\inn{\omega^{-1/2}\phi}{h^{\omega^{-1}}_K}_{\omega^{-1}}C_{K-}(\omega)\inn{h^\omega_{K-}}{g}_\omega\tag{A11}\\
    +&\sum_KD_K(\omega^{-1})\inner{\phi\omega^{-1/2}}_KC_{K-}(\omega)\inn{h^\omega_{K-}}{g}_\omega\tag{A12}
\end{align*}
and
\begin{align*}
    &\sum_K\innerHaar{\omega^{\frac{-1}{2}}\phi}{K}D_{K-}(\omega)\inn{h^1_{K-}}{g}_\omega\tag{A2}\\
    =&C_K(\omega^{-1})\inn{\omega^{-1/2}\phi}{h^{\omega^{-1}}_K}_{\omega^{-1}}D_{K-}(\omega)\inn{h^1_{K-}}{g}_\omega\tag{A21}\\
    +&D_K(\omega^{-1})\inner{\phi\omega^{-1/2}}_KD_{K-}(\omega)\inn{h^1_{K-}}{g}_\omega.\tag{A22}
\end{align*}
\paragraph{For A11}
\begin{align*}
    &\sum_KC_K(\omega^{-1})\inn{\omega^{-1/2}\phi}{h^{\omega^{-1}}_K}_{\omega^{-1}}C_{K-}(\omega)\inn{h^\omega_{K-}}{g}_\omega\\
    \lesssim&\sum_K\sqrt{\inner{\omega^{-1}}_K\inner{\omega}_K}\inn{\omega^{-1/2}\phi}{h^{\omega^{-1}}_K}_{\omega^{-1}}\inn{h^\omega_{K-}}{g}_\omega\quad\textit{by (\ref{Cestimate}), (\ref{C-estimate})}\\
    \leq&\bracket{\sum_K\sqrt{\inner{\omega^{-1}}_K\inner{\omega}_K}\inn{\omega^{-1/2}\phi}{h^{\omega^{-1}}_K}^2_{\omega^{-1}}}^{1/2}\bracket{\sum_K\sqrt{\inner{\omega^{-1}}_K\inner{\omega}_K} \inn{h^\omega_{K-}}{g}^2_\omega}^{1/2}\\
    \leq& \Atwo^{1/2}\norm{\omega^{1/2}\phi}_{\omega^{-1}}\norm{g}_\omega\\
    =&\Atwo^{1/2}\norm{\phi}\norm{g}_\omega.
\end{align*}
\paragraph{For A12}
By Cauchy-Schwarz, we have
\begin{align*}
    &\sum_KD_K(\omega^{-1})\inner{\phi\omega^{-1/2}}_KC_{K-}(\omega)\inn{h^\omega_{K-}}{g}_\omega\\
    \leq& \bracket{\sum_KD_K(\omega^{-1})^2\inner{\phi\omega^{-1/2}}^2_KC^2_{K-}(\omega)}^{1/2}\norm{g}_\omega.\\
\end{align*}
By weighted Carleson embedding(\ref{wCarleson}) applied to $\phi$ we have that
\begin{align*}
    & \bracket{\sum_KD_K(\omega^{-1})^2\inner{\phi\omega^{-1/2}}^2_KC^2_{K-}(\omega)}\\
    \lesssim &\bracket{\sum_K \frac{\innerHaar{\omega^{-1}}{K}^2}{\inner{\omega^{-1}}^2_K}\inner{\omega}_K\inner{\phi\omega^{-1/2}}^2_K}\\
    \leq& \norm{\phi}^2\sup_J\frac{1}{|J|\inner{\omega^{-1}}_J}\sum_{K\subset J}\inner{\omega^{-1}}^2_K\frac{\innerHaar{\omega^{-1}}{K}^2}{\inner{\omega^{-1}}^2_K}\inner{\omega}_K\\
    =&\norm{\phi}^2\sup_J\frac{1}{|J|\inner{\omega^{-1}}_J}\sum_{K\subset J}\innerHaar{\omega^{-1}}{K}^2\inner{\omega}_K\\
    \leq& \norm{\phi}^2\Atwo^2\quad\textit{by (\ref{Swestimate})}.
\end{align*}
Thus
\begin{equation*}
    \sum_KD_K(\omega^{-1})\inner{\phi\omega^{-1/2}}_KC_{K-}(\omega)\inn{h^\omega_{K-}}{g}_\omega\lesssim \Atwo \norm{\phi}\norm{g}_\omega.
\end{equation*}
\paragraph{For A21}
By Cauchy-Schwarz, we have
\begin{align*}
    &\sum_KC_K(\omega^{-1})\inn{\omega^{-1/2}\phi}{h^{\omega^{-1}}_K}D_{K-}(\omega)\inn{h^1_{K-}}{g}_\omega\\
    \leq &\norm{\omega^{1/2}\phi}_{\omega^{-1}}\bracket{\sum_KC_K(\omega^{-1})^2D_{K-}(\omega)^2\inn{h^1_{K-}}{g}^2_\omega}^{1/2}\\
    =&\norm{\phi}\bracket{\sum_KC_K(\omega^{-1})^2D_{K-}(\omega)^2\inn{h^1_{K-}}{g}^2_\omega}^{1/2}\\
    \lesssim&\norm{\phi}\bracket{\sum_K\inner{\omega^{-1}}_K\frac{\innerHaar{\omega}{K-}^2}{\inner{\omega}^2_{K-}}\inner{g\omega}^2_{K-}}^{1/2}\quad\textit{by (\ref{Cestimate}).}
\end{align*}
By weighted Carleson embedding(\ref{Carlesonem}) applied to $g$ we have that
\begin{align*}
    &\sum_K\inner{\omega^{-1}}_K\frac{\innerHaar{\omega}{K-}^2}{\inner{\omega}^2_{K-}}\inner{g\omega}^2_{K-}\\
    \leq &\norm{g}^2_\omega\sup_J\frac{1}{|J|\inner{\omega}_J}\sum_{K\subset J}\inner{\omega}^2_{K-}\inner{\omega^{-1}}_K\frac{\innerHaar{\omega}{K-}^2}{\inner{\omega}^2_{K-}}\\
    =&\norm{g}^2_\omega\sup_J\frac{1}{|J|\inner{\omega}_J}\sum_{K\subset J}\inner{\omega^{-1}}_K\innerHaar{\omega}{K-}^2\\
    \lesssim&\norm{g}^2_\omega\Atwo^2.
\end{align*}
The last inequality comes from the modified square function (Theorem \ref{modified}).
Hence 
\begin{equation*}
    \sum_KC_K(\omega^{-1})\inn{\omega^{-1/2}\phi}{h^{\omega^{-1}}_K}D_{K-}(\omega)\inn{h^1_{K-}}{g}_\omega\lesssim \Atwo\norm{g}_\omega\norm{\phi}.
\end{equation*}
\paragraph{For A22}
This term we apply bilinear Carleson embedding (\ref{Bilinear Imbedding Theorem}) and (\ref{wineq}), (\ref{iwineq}), and (\ref{ineq}). 
\begin{align*}
    &D_K(\omega^{-1})\inner{\phi\omega^{-1/2}}_KD_{K-}(\omega)\inn{h^1_{K-}}{g}_\omega\\
    \leq&\frac{|\innerHaar{\omega^{-1}}{K}\innerHaar{\omega}{K-}|}{\inner{\omega^{-1}_K}\inner{\omega}_{K-}}\inn{\phi\omega^{1/2}}{h^1_K}_{\omega^{-1}}\inn{h^1_{K-}}{g}_\omega\\
    \lesssim &\Atwo\norm{\phi\omega^{1/2}}_{\omega^{-1}}\norm{g}_\omega\\
    = &\Atwo\norm{\phi}\norm{g}_\omega.
\end{align*}
Thus all the terms are bounded by $\Atwo\norm{\phi}\norm{g}_\omega$. Hence the proof is complete.
\section{Preliminary for Haar multiplier}\label{PreMul}
Let $L^2\equiv L^2(\mathbb{R}^d)$ be the space of square integrable functions over $\mathbb{R}^d$. A weight is a nonnegative local integrable function i.e. $\omega\in L^1_{loc}(\mathbb{R}^d)$ and $$\omega:\mathbb{R}^d\rightarrow \mathbb{R}_{\geq 0}.$$
Also $\inn{}{}:=\inn{}{}_{L^2(\mathbb{R}^d)}$, and $\inn{}{}_\omega:=\inn{}{}_{L^2(\mathbb{R}^d,\omega)}$; $\norm{\cdot}=\norm{\cdot}_{L^2(\mathbb{R}^d)}$, and $\norm{\cdot}_\omega=\norm{\cdot}_{L^2(\mathbb{R}^d,\omega)}$. Let $\mathcal{D}$ denote the dyadic grid in $\mathbb{R}^d$ and 
\begin{equation*}
    \Gamma_d\equiv\{0,1\}^d\backslash\{(1,\dots,1)\}.
\end{equation*}

Now we recall the Wilson system \cite{wilson}
\subsection{Wilson system}
Define $\mathcal{C}_1(I):\{J\in\mathcal{D} | J\subset I\}$. The following lemma is used to construct our Haar system and Carleson embedding in the later.
\begin{lem}
(Wilson, \cite{wilson})
Let $I\in\mathcal{D}$. Then there are $2^d-1$ pairs of sets $\{(E^1_{\alpha,I},E^2_{\alpha,I})\}$ such that
\begin{description}
\item[1]For each $\alpha\in\Gamma_d$, $|E^1_{\alpha,I}|=|E^2_{\alpha,I}|$;
\item[2]For each $\alpha$ and $s=1,2$, $E^s_{\alpha,I}$ is a non-empty union of cubes from $\mathcal{C}_1(I)$;
\item[3]For each $\alpha$, $E^1_{\alpha,I}\cap E^2_{\alpha,I}\neq\emptyset$;
\item[4]For each $\alpha\neq\beta$ one of the following must hold:\\
\quad(a)$E^1_{\alpha,I}\cup  E^2_{\alpha,I}$ is entirely contained in either $E^1_{\beta,I}$ or $E^2_{\beta,I}$;\\
\quad(b)$E^1_{\beta,I}\cup  E^2_{\beta,I}$ is entirely contained in either $E^1_{\alpha,I}$ or $E^2_{\alpha,I}$;\\
\quad(c)$\bracket{E^1_{\beta,I}\cup  E^2_{\beta,I}}\cap\bracket{E^1_{\beta,I}\cup  E^2_{\beta,I}}=\emptyset$.
\end{description}
\end{lem}
Let $E_{\alpha,I}=E^1_{\beta,I}\cup  E^2_{\beta,I}$. Now we introduce Wilson's Haar system of $L^2(\omega)$. For $\alpha\in\Gamma_d,I\in\mathcal{D}$, define
\begin{equation*}
    h^{\omega,\alpha}_I\equiv\frac{1}{\sqrt{\omega(E_{\alpha,I})}}\bracket{\frac{\sqrt{\omega(E^1_{\alpha,I})}}{\sqrt{\omega(E^2_{\alpha,I})}}1_{E^2_{\alpha,I}}-\frac{\sqrt{\omega(E^2_{\alpha,I})}}{\sqrt{\omega(E^1_{\alpha,I})}}1_{E^1_{\alpha,I}}}.
\end{equation*}
For a fixed $J\in\mathcal{D}$, set
\begin{align*}
    \innerHaar{f}{J,\alpha}\equiv\inn{f}{h^\alpha_J}\quad\forall\alpha\in\Gamma_d,\\
    h_J\equiv\bracket{\innerHaar{f}{J,\alpha}}_{\alpha\in\Gamma_d}.
\end{align*}
For a giving set $E$, define $h_E^1\equiv\frac{1}{|E|}1_E$, and $\inner{f}_E=\inn{f}{h^1_E}$.
Then we define Haar multiplier $T_\sigma$. Let $\{\sigma_{I,\alpha}\}_{\alpha\in\Gamma_d,I\in\mathcal{D}}$ denote a number of sequence with $\sigma_{I,\alpha}\in\{\pm 1\}$. The Haar multiplier $T_\sigma$ is defined by
\begin{equation*}
    T_\sigma\equiv\Asum \sigma_{I\alpha}\innerHaar{f}{I,\alpha}h^\alpha_I.
\end{equation*}
A fundamental property for the Wilson Haar system is that, for $f,g\in L^2$
\begin{equation*}
    \innerHaar{fg}{J,\beta}=\Asum\innerHaar{f}{I,\alpha}\innerHaar{g}{I,\alpha}\inn{h^1_{E_{\alpha,I}}}{h^\beta_J}+\innerHaar{f}{J,\beta}\inner{g}_{E_\beta,J}+\innerHaar{g}{J,\beta}\inner{f}_{E_{\beta,J}}.
\end{equation*}

Now we consider the following dyadic operators \cite{Multiplier}. Giving a sequence of numbers $a=\{a_{I,\alpha}\}_{I\in\mathcal{D},\alpha\in\Gamma_d}$ indexed by $I\in\mathcal{D}$ and $\alpha\in\Gamma_d$, we have the following paraproduct type operators
\begin{align*}
    &P^{(0,0)}_af\equiv\Asum a_{I\,\alpha}\innerHaar{f}{I,\alpha}h^\alpha_I\\
    &P^{(0,1)}_af\equiv\Asum a_{I\,\alpha}\inn{f}_{E_{\alpha,I}}h^\alpha_I\\
    &P^{(1,0)}_af\equiv\Asum a_{I\,\alpha}\innerHaar{f}{I,\alpha}h^1_{E_{\alpha,I}}.
\end{align*}
It is not hard to see that the operator $M_g$ of multiplication by $g$ can formally be written as
\begin{equation*}
    M_gf=P^{(0,0)}_{\inner{g}}f+P^{(0,1)}_{\widehat{g}}f+P^{(1,0)}_{\widehat{g}}f,
\end{equation*}
where $\inner{g}\equiv\{\inner{g}_{E_{\alpha,I}}\}_{I\in\mathcal{D},\alpha\in\Gamma_d}$ and
$\widehat{g}\equiv\{\widehat{g}_{E_{\alpha,I}}\}_{I\in\mathcal{D},\alpha\in\Gamma_d}$. 

Finally we have disbalanced Haar functions: Fixing a dyadic cube $J$, a weight $\omega$ on $\mathbb{R}^d$, and $\beta\in\Gamma$, we set
\begin{equation*}
    C_J(\omega,\beta)\equiv\sqrt{\frac{\inner{\omega}_{E^1_{\beta,J}}\inner{\omega}_{E^2_{\beta,J}}}{\inner{\omega}_{E_{\beta,J}}}}\quad\textit{ and }\quad D_J(\omega,\beta)\equiv\frac{\innerHaar{\omega}{J,\beta}}{\inn{\omega}_{E_{\beta,J}}}.
\end{equation*}
Then we have
\begin{equation*}
    h^\beta_J=C_J(\omega,\beta)h^{\omega,\beta}_J+D_J(\omega,\beta)h^1_{E_{\beta,J}}.
\end{equation*}
Also we have estimate of $C_J(\omega,\beta)$
\begin{equation*}
    C_J(\omega,\beta)^2\leq 4\inner{\omega}_{E_{\beta,J}}\leq 2^{d+1}\inner{\omega}_J.
\end{equation*}
\subsection{Carleson embedding theorem}
In theorem 4.3 \cite{Chung}, we have following result
\begin{thm}\label{Carlemndim}
(Modified Carleson Embedding Theorem). Let $\omega$ be a weight on $\mathbb{R}^d$ and let $\{a_{\alpha,I}\}_{I\in\mathcal{D},\alpha\in\Gamma_d}$ be a sequence of nonnegative numbers. Then, there is a constant $A>0$ such that
\begin{equation*}
    \frac{1}{|E_{\alpha,I}|}\sum_{J\subset I}\sum_{\beta:E_{\beta,J}\subset E_{\alpha,I}}a_{\beta,J}\inner{\omega}^2_\Eb\leq A\inner{\omega}_{E_{\alpha,I}}\quad \forall I\in\mathcal{D},\alpha\in\Gamma_d,
\end{equation*}
if and only if
\begin{equation*}
    \Asum a_{\alpha,I}\inner{\omega^{\frac{1}{2}}f}_{E_{\alpha,I}}\lesssim A\norm{f}^2.
\end{equation*}
\end{thm}
Besides, we have bilinear Carleson embedding from \cite{Chung}
\begin{thm}\label{Multicarl}
(Multivariable Verson of the Bilinear Carleson Embedding Theorem) Let $\omega, \nu$ be weights so that $\inner{\omega}_{E_{\alpha,I}}\inner{\omega^{-1}}_{E_{\alpha,I}}\leq A$, for all $\alpha\in\Gamma_d$ and $I\in\mathcal{D}$, and $\{a_{I,\alpha}\}_{I,\alpha}$ is a sequence of nonnegative numbers such that the three inequalities below hold with some constant $A>0$,
\begin{align*}
    &\Carlsum a_{I,\alpha}\inner{\nu}_{\Ea}\leq A\inner{\nu}_{\Ea}\\
    &\Carlsum a_{I,\alpha}\inner{\omega}_{\Ea}\leq A\inner{\omega}_{\Ea}\\
    &\Carlsum a_{I,\alpha}\inner{\omega}_{\Ea}\inner{\nu}_{\Ea}\leq A.
\end{align*}
Then for all $f\in L^2(\omega^{-1})$, $g\in L^2(\nu^{-1})$
\begin{equation*}
    \Asum a_{I,\alpha}\inner{f}_{\Ea}\inner{g}_{\Ea}\leq CA\norm{f}_{\omega^{-1}}\norm{g}_{\nu^{-1}}
\end{equation*}
holds with some constant $C>0$.
\end{thm}
Also we have the following inequalities from (6.3),(6.4) \cite{Chung}
\begin{align}
    &\Carlsum|\innerHaar{\omega}{J,\beta}\innerHaar{\omega^{-1}}{J,\beta}|\leq C(n)\Atwo\label{ninequal}\\
    &\wCarlsum\frac{|\innerHaar{\omega}{J,\beta}\innerHaar{\omega^{-1}}{J,\beta}|}{\inner{\omega^{-1}}_{E_{\beta,J}}}\leq C(n)\Atwo\label{nwinequal}\\
    &\wCarlsumne\frac{|\innerHaar{\omega}{J,\beta}\innerHaar{\omega^{-1}}{J,\beta}|}{\inner{\omega}_{E_{\beta,J}}}\leq C(n)\Atwo.\label{winequal}
\end{align}
\subsection{Square function}
\begin{thm}\label{Square}
\cite{Multiplier}For $\omega\in A_2$,
$$\norm{Sf}_\omega\lesssim\Atwo\norm{f}_\omega.$$
Using this equality to $f=\omega^\frac{-1}{2}1_{E_{\alpha,I}}$ for $I\in\mathcal{D},\alpha\in\Gamma_d$ yields the following:
$$\sum_{J\subseteq I}\sum_{\beta:E_{\beta,J}\subseteq E_{\alpha,I}}|\innerHaar{\omega^\frac{-1}{2}}{J,\beta}|^2\inner{\omega}_J\lesssim\Atwo^2|E_{\alpha,I}|\quad\forall I\in\mathcal{D},\alpha\in\Gamma_d.$$
\end{thm}
\section{Main result for Haar Multiplier}\label{Haarmulti}
What we will prove is the following 

$$\norm{P^{(0,1)}_{\widehat{\omega^{\frac{1}{2}}}}T_\sigma P^{(1,0)}_{\widehat{\omega^{\frac{-1}{2}}}}\phi}\lesssim \Atwo\norm{\phi},$$

and it is equivalent to show

$$\norm{T_\sigma P^{(1,0)}_{\widehat{\omega^{\frac{-1}{2}}}}\phi}_\omega\lesssim\Atwo\norm{\phi}$$
since we use similar idea as for the Haar shift above and the two terms are already proved in the work of K. Bickel, E. Sawyer and B. Wick\cite{Multiplier}.
Following the definition, we see
\begin{align*}
  P^{(1,0)}_{\widehat{\omega^{\frac{-1}{2}}}}\phi=& \Asum\innerHaar{\omega^{\frac{-1}{2}}}{I,\alpha}\innerHaar{\phi}{I,\alpha}\wHaara^1,\\
  T_\sigma P^{(1,0)}_{\widehat{\omega^{\frac{-1}{2}}}}\phi=&\Asum\innerHaar{\omega^{\frac{-1}{2}}}{I,\alpha}\innerHaar{\phi}{I,\alpha}\Bsum\inn{\wHaara^1}{h^\beta_J}h^\beta_J\cdot\sigma_{J,\beta}\\
  =& \Bsum\cbracket{\sum_{I\subset J}\sum_{E_{\alpha,I}\subsetneq E_{\beta,J}}\innerHaar{\omega^{\frac{-1}{2}}}{I,\alpha}\innerHaar{\phi}{I,\alpha}\inn{h^1_{E_{\alpha,I}}}{h^\beta_J}}h^\beta_J\cdot\sigma_{J,\beta}\\
  =&\Bsum\cbracket{\innerHaar{\omega^{\frac{-1}{2}}\phi}{J,\beta}-\innerHaar{\omega^{\frac{-1}{2}}}{J,\beta}\inner{\phi}_\Eb-\innerHaar{\phi}{J,\beta}\inner{\omega^{\frac{-1}{2}}}_\Eb}h^\beta_J\cdot\sigma_{J,\beta}.
\end{align*}
Now we turn to further decomposition
\begin{prop}\label{whaardec}
$\inn{h^\beta_J}{\omega g }=C_J(\omega,\beta)\inn{g}{h^{\omega,\beta}_J}_\omega+D_J(\omega,\beta)\inner{\omega g}_\Eb$
\end{prop}
Using proposition (\ref{whaardec}) and duality, given a $g\in L^2(\omega)$ 
\begin{align*}
&\inn{T_\sigma P^{(1,0)}_{\widehat{\omega^{\frac{-1}{2}}}}\phi}{g}_\omega\\
=&\Bsum\cbracket{\innerHaar{\omega^{\frac{-1}{2}}\phi}{J,\beta}-\innerHaar{\omega^{\frac{-1}{2}}}{J,\beta}\inner{\phi}_\Eb-\innerHaar{\phi}{J,\beta}\inner{\omega^{\frac{-1}{2}}}_\Eb}\sigma_{J,\beta}\cdot\inn{h^\beta_J}{\omega g}\\
=&\Bsum\cbracket{\innerHaar{\omega^{\frac{-1}{2}}\phi}{J,\beta}-\innerHaar{\omega^{\frac{-1}{2}}}{J,\beta}\inner{\phi}_\Eb-\innerHaar{\phi}{J,\beta}\inner{\omega^{\frac{-1}{2}}}_\Eb}\sigma_{J,\beta}\\
\cdot& \cbracket{C_J(\omega,\beta)\inn{g}{h^{\omega,\beta}_J}_\omega+D_J(\omega,\beta)\inner{\omega g}_\Eb}.
\end{align*}
We split it into the following 6 terms i.e.
\begin{align*}
    &\sigma_{J,\beta}\Bsum\innerHaar{\omega^{\frac{-1}{2}}\phi}{J,\beta}C_J(\omega,\beta)\inn{g}{h^{\omega,\beta}_J}_\omega \tag{A1}\\
    &\sigma_{J,\beta}\Bsum\innerHaar{\phi}{J,\beta}\inner{\omega^{\frac{-1}{2}}}_\Eb C_J(\omega,\beta)\inn{g}{h^{\omega,\beta}_J}_\omega\tag{B1}\\
    &\sigma_{J,\beta}\Bsum\innerHaar{\omega^{\frac{-1}{2}}}{J,\beta}\inner{\phi}_\Eb C_J(\omega,\beta)\inn{g}{h^{\omega,\beta}_J}_\omega\tag{C1}\\
    &\sigma_{J,\beta}\Bsum\innerHaar{\omega^{\frac{-1}{2}}\phi}{J,\beta} D_J(\omega,\beta)\inner{\omega g}_\Eb\tag{A2}\\
    &\sigma_{J,\beta}\Bsum\innerHaar{\phi}{J,\beta}\inner{\omega^{\frac{-1}{2}}}_\Eb D_J(\omega,\beta)\inner{\omega g}_\Eb\tag{B2}\\
    &\sigma_{J,\beta}\Bsum\innerHaar{\omega^{\frac{-1}{2}}}{J,\beta}\inner{\phi}_\Eb D_J(\omega,\beta)\inner{\omega g}_\Eb\tag{C2}.\\
\end{align*}
Since $|\sigma_{I,\alpha}|=1$, it suffices to show all the 6 terms (A1, A2,B1, B2, C1, C2) are bounded by $\Atwo\norm{\phi}\norm{g}_\omega.$ 
\paragraph{For B1} Applying Cauchy-Schwarz twice
\begin{align*}
    &\Bsum\innerHaar{\phi}{J,\beta}\inner{\omega^{\frac{-1}{2}}}_\Eb C_J(\omega,\beta)\inn{g}{h^{\omega,\beta}_J}_\omega\\
    \lesssim&\norm{\phi}\bracket{\Bsum\inner{\omega^{\frac{-1}{2}}}^2_\Eb C_J(\omega,\beta)^2\inn{g}{h^{\omega,\beta}_J}^2_\omega}^\frac{1}{2}\\
    \lesssim&\norm{\phi}\bracket{\Bsum\inner{\omega^{-1}}_\Eb\inner{\omega}_\Eb\inn{g}{h^{\omega,\beta}_J}^2_\omega}^\frac{1}{2}\\
    \lesssim &\Atwo\norm{\phi}\norm{g}_\omega.
\end{align*}
\paragraph{For C1} By Cauchy-schwarz inequality and Carleson embedding theorem (\ref{Carlemndim})
\begin{align*}
    &\Bsum \innerHaar{\omega^{\frac{-1}{2}}}{J,\beta}\inner{\phi}_\Eb C_J(\omega,\beta)\inn{g}{h^{\omega,\beta}_J}_\omega\\
    \leq&\norm{g}_\omega\bracket{\Bsum\innerHaar{\omega^{\frac{-1}{2}}}{J,\beta}^2\inner{\phi}^2_\Eb C_J(\omega,\beta)^2}^\frac{1}{2}\\
    \lesssim&\norm{g}_\omega\bracket{\Bsum\innerHaar{\omega^\frac{-1}{2}}{J,\beta}^2\inner{\omega}_J\inner{\phi}^2_\Eb}^\frac{1}{2}\\
    \lesssim & \norm{\phi}\norm{g}_\omega\bracket{\sup_I\frac{1}{|I|}\Bsum \innerHaar{\omega^\frac{-1}{2}}{J,\beta}^2\inner{\omega}_J}^\frac{1}{2}\\
    \lesssim&\Atwo \norm{\phi}\norm{g}_\omega.
\end{align*}
The last inequality comes from estimation deduced from the square function (\ref{Square}).
\paragraph{For B2}
By Carleson embedding (Theorem \ref{Carlemndim}), we have
\begin{align*}
    &\Bsum\innerHaar{\phi}{J,\beta}\inner{\omega^{\frac{-1}{2}}}_\Eb D_J(\omega,\beta)\inner{\omega g}_\Eb\\
    \leq&\norm{\phi}\bracket{\Bsum\inner{\omega^{\frac{-1}{2}}}_{\Eb}^2D_J(\omega,\beta)^2\inner{\omega g}^2_\Eb}^\frac{1}{2}\quad\textit{by Cauchy-schwarz}\\
    \leq&\norm{\phi}\norm{g}_\omega\bracket{\wCarlsum\inner{\omega^{-1}}_\Eb\innerHaar{\omega}{J,\beta}^2}^{\frac{1}{2}}\\
    \lesssim& \Atwo\norm{\phi}\norm{g}_\omega.
\end{align*}
The last inequality comes from square function estimate (Theorem \ref{Square}).
\paragraph{For C2}
In this part, we need to split it into two cases. Since we have proposition \ref{whaardec}.
Thus
\begin{equation*}
    \innerHaar{\omega^{\frac{-1}{2}}}{J,\beta}=C_J(\omega^{-1},\beta)\inn{\omega^\frac{1}{2}}{h^{\omega^{-1},\beta}_J}_{\omega^{-1}}+D_J(\omega^{-1},\beta)\inner{\omega^\frac{-1}{2}}_\Eb.
\end{equation*}
\begin{align*}
    &\Bsum\innerHaar{\omega^{\frac{-1}{2}}}{J,\beta}\inner{\phi}_\Eb D_J(\omega,\beta)\inner{\omega g}_\Eb\\
    =&\Bsum C_J(\omega^{-1},\beta)\inn{\omega^\frac{1}{2}}{h^{\omega^{-1},\beta}_J}_{\omega^{-1}}\inner{\phi}_\Eb D_J(\omega,\beta)\inner{\omega g}_\Eb\tag{C21}\\
    +&\Bsum D_J(\omega^{-1},\beta)\inner{\omega^\frac{-1}{2}}_\Eb \inner{\phi}_\Eb D_J(\omega,\beta)\inner{\omega g}_\Eb.\tag{C22}
\end{align*}
\paragraph{For C21}
By Cauchy-Schwarz, we have
\begin{align*}
    &\Bsum C_J(\omega^{-1},\beta)\inn{\omega^\frac{1}{2}}{h^{\omega^{-1},\beta}_J}_{\omega^{-1}}\inner{\phi}_\Eb D_J(\omega,\beta)\inner{\omega g}_\Eb\\
    \leq&\bracket{\Bsum\inn{\omega^\frac{1}{2}}{h^{\omega^{-1},\beta}_J}_{\omega^{-1}}^2\inner{\phi}_\Eb^2 }^\frac{1}{2}\bracket{\Bsum\inner{\omega^{-1}}_\Eb D_J(\omega,\beta)^2\inner{\omega g}^2_\Eb}^\frac{1}{2}.
\end{align*}
By Carleson embedding,
\begin{align*}
    &\bracket{\Bsum\inn{\omega^\frac{1}{2}}{h^{\omega^{-1},\beta}_J}_{\omega^{-1}}^2\inner{\phi}_\Eb^2 }\\
    \lesssim &\norm{\phi}^2\Carlsum \inn{\omega^\frac{1}{2}}{h^{\omega^{-1},\beta}_J}_{\omega^{-1}}^2\\
    \leq& \norm{\phi}^2
\end{align*}
By Modified Carleson embedding theorem \ref{Carlemndim}, 
\begin{align*}
    &\bracket{\Bsum\inner{\omega^{-1}}_\Eb D_J(\omega,\beta)^2\inner{\omega g}^2_\Eb}\\
    \lesssim &\norm{g}_\omega\wCarlsum\inner{\omega^{-1}}_\Eb\innerHaar{\omega}{J,\beta}^2\\
    \lesssim &\Atwo\norm{g}_\omega.
\end{align*}
The last inequality comes from square function estimate.
\paragraph{For C22}
We apply Carleson embedding for $\phi$ and weight Carleson embedding (\ref{Carlemndim}) for $g$.
\begin{align*}
    &\Bsum D_J(\omega^{-1},\beta)\inner{\omega^\frac{-1}{2}}_\Eb \inner{\phi}_\Eb D_J(\omega,\beta)\inner{\omega g}_\Eb\\
    \leq&\bracket{\Bsum\innerHaar{\omega^{-1}}{J,\beta}\innerHaar{\omega}{J,\beta}\inner{\phi}_\Eb}^\frac{1}{2}\bracket{\Bsum\frac{\innerHaar{\omega^{-1}}{J,\beta}\innerHaar{\omega}{J,\beta}}{\inner{\omega^{-1}}_\Eb\inner{\omega}_\Eb^2}\inner{\omega g}^2_\Eb}^\frac{1}{2}\quad\textit{by Cauchy-schwarz}\\
    \lesssim&\norm{\phi}\bracket{\Carlsum|\innerHaar{\omega}{J,\beta}\innerHaar{\omega^{-1}}{J,\beta}|}^\frac{1}{2}\norm{g}_\omega\\
    &\bracket{\wCarlsum\frac{|\innerHaar{\omega}{J,\beta}\innerHaar{\omega^{-1}}{J,\beta}|}{\inner{\omega^{-1}}_{E_{\beta,J}}}}^\frac{1}{2}\\
    \lesssim&\Atwo\norm{\phi}\norm{g}_\omega\quad\textit{ by (\ref{ninequal}) and (\ref{nwinequal}).}
\end{align*}
\paragraph{For part A:}
Again, we decompose $A1, A2$ into $A11, A12, A21, A22$. Note
\begin{equation*}
    \innerHaar{\omega^\frac{-1}{2}}{J,\beta}=C_J(\omega^{-1},\beta)\inn{\omega^\frac{1}{2}\phi}{h^{\omega^{-1},\beta}_J}_{\omega^{-1}}+D_J(\omega^{-1},\beta)\inner{\omega^\frac{-1}{2}\phi}_\Eb.
\end{equation*}
Hence
\begin{align*}
    &\Bsum\innerHaar{\omega^{\frac{-1}{2}}\phi}{J,\beta}C_J(\omega,\beta)\inn{g}{h^{\omega,\beta}_J}_\omega\tag{A1}\\
    =&\Bsum C_J(\omega^{-1},\beta)\inn{\omega^\frac{1}{2}\phi}{h^{\omega^{-1},\beta}_J}_{\omega^{-1}}C_J(\omega,\beta)\inn{g}{h^{\omega,\beta}_J}_\omega\tag{A11}\\
    +&\Bsum D_J(\omega^{-1},\beta)\inner{\omega^\frac{-1}{2}\phi}_\Eb C_J(\omega,\beta)\inn{g}{h^{\omega,\beta}_J}_\omega\tag{A12}\\
\end{align*}
and 
\begin{align*}
    &\Bsum\innerHaar{\omega^{\frac{-1}{2}}\phi}{J,\beta} D_J(\omega,\beta)\inner{\omega g}_\Eb\tag{A2}\\
    =&\Bsum C_J(\omega^{-1},\beta)\inn{\omega^\frac{1}{2}\phi}{h^{\omega^{-1},\beta}_J}_{\omega^{-1}}D_J(\omega,\beta)\inner{\omega g}_\Eb\tag{A21}\\
    +&\Bsum D_J(\omega^{-1},\beta)\inner{\omega^\frac{-1}{2}\phi}_\Eb D_J(\omega,\beta)\inner{\omega g}_\Eb.\tag{A22}
\end{align*}
\paragraph{For A11}
By Cauchy-Schwarz:
\begin{align*}
    &\Bsum C_J(\omega^{-1},\beta)\inn{\omega^\frac{1}{2}\phi}{h^{\omega^{-1},\beta}_J}_{\omega^{-1}}C_J(\omega,\beta)\inn{g}{h^{\omega,\beta}_J}_\omega\\
    \lesssim&\bracket{\Bsum\sqrt{\inner{\omega^{-1}}_{E_{\beta,J}}\inner{\omega}_{E_{\beta,J}}}\inn{\omega^\frac{1}{2}\phi}{h^{\omega^{-1},\beta}_J}_{\omega^{-1}}^2}^\frac{1}{2}\bracket{\Bsum\sqrt{\inner{\omega^{-1}}_{E_{\beta,J}}\inner{\omega}_{E_{\beta,J}}}\inn{g}{h^{\omega,\beta}_J}_\omega^2}^\frac{1}{2}\\
    \leq& \Atwo^\frac{1}{2}\norm{\phi}\norm{g}_{\omega}.
\end{align*}
\paragraph{For A12}
By Cauchy-Schwarz:
\begin{align*}
    &\Bsum D_J(\omega^{-1},\beta)\inner{\omega^\frac{-1}{2}\phi}_\Eb C_J(\omega,\beta)\inn{g}{h^{\omega,\beta}_J}_\omega\\
    \leq &\norm{g}_\omega\bracket{D_J(\omega^{-1},\beta)^2\inner{\omega^\frac{-1}{2}\phi}^2_\Eb C_J(\omega,\beta)^2}^\frac{1}{2}\\
    \leq &\norm{g}_\omega\norm{\phi}\bracket{\wCarlsumne \inner{\omega}_\Eb \innerHaar{\omega^{-1}}{J,\beta}^2}^\frac{1}{2}\quad\textit{by \ref{Carlemndim}}\\
    \lesssim& \Atwo \norm{g}_\omega\norm{\phi}.
\end{align*}
The last inequality comes from square function estimate.
\paragraph{For A21}
By Cauchy-Schwarz:
\begin{align*}
    &\Bsum C_J(\omega^{-1},\beta)\inn{\omega^\frac{1}{2}\phi}{h^{\omega^{-1},\beta}_J}_{\omega^{-1}}D_J(\omega,\beta)\inner{\omega g}_\Eb\\
    \leq&\norm{\omega^\frac{1}{2}\phi}_{\omega^{-1}}\bracket{\Bsum C_J(\omega^{-1},\beta)^2D_J(\omega,\beta)^2\inner{\omega g}^2_\Eb}^\frac{1}{2}\\
    \lesssim&\norm{\phi}\norm{g}_\omega\bracket{\wCarlsum \inner{\omega^{-1}_\Eb}^2\innerHaar{\omega}{J,\beta}}^\frac{1}{2}\\
    \lesssim&\Atwo\norm{\phi}\norm{g}_\omega.\quad\textit{by Theorem \ref{Square}}.
\end{align*}
\paragraph{For A22}
By Cauchy-Schwarz and Multivariable version of the bilinear Carleson embedding (Theorem \ref{Multicarl}) with inequalities (\ref{ninequal}) (\ref{wineq}) (\ref{nwinequal}), we have
\begin{align*}
    &\Bsum D_J(\omega^{-1},\beta)\inner{\omega^\frac{-1}{2}\phi}_\Eb D_J(\omega,\beta)\inner{\omega g}_\Eb\\
    =&\bracket{\frac{\innerHaar{\omega^{-1}}{J,\beta}\innerHaar{\omega}{J,\beta}}{\inner{\omega }_\Eb\inner{\omega^{-1} }_\Eb}\inner{\omega^\frac{-1}{2}\phi}_\Eb\inner{\omega g}_\Eb}\\
    \leq &\Atwo\norm{\omega^\frac{-1}{2}\phi}_{\omega}\norm{\omega g}_{\omega^{-1}}\\
    =&\Atwo\norm{\phi}\norm{g}_\omega.
\end{align*}
Since all the terms are bounded by $\Atwo\norm{\phi}\norm{g}_\omega$. Hence the proof is complete.

\end{document}